\documentstyle[11pt,epsfig]{article}
\begin{document}

\setlength{\arraycolsep}{2pt}
\newcounter{itemnbr}
\newcommand\eq[1] {(\ref{#1})}
\newcommand\eqp[1] {(\ref{eq:#1})}
\newcommand\defn[1] {\ref{def:#1}}
\newcommand\thm[1] {\ref{thm:#1}}
\newcommand\sect[1] {\ref{sec:#1}}
\newcommand\fig[1] {\ref{fig:#1}}
\newcommand\labeq[1] {\label{eq:#1}}
\newcommand\labdefn[1] {\label{def:#1}}
\newcommand\labthm[1] {\label{thm:#1}}
\newcommand\labsect[1] {\label{sec:#1}}
\newcommand\labfig[1] {\label{fig:#1}}
\newcommand\etal{{\it et al.}}
\newtheorem{definition}{Definition}[section]
\newtheorem{lemma}{Lemma}[section]
\newtheorem{theorem}{Theorem}[section]
\newtheorem{corollary}{Corollary}[section]
\newcommand\proof{{\it Proof:}\quad}
\newcommand\remark{{\it Remark:}\quad}
\newcommand\qed{\ \rule[-0.2ex]{0.3em}{1.5ex}}
\newcommand{\bfm}[1]{\mbox{\boldmath ${#1}$}}
\newcommand{\nonum}{\nonumber \\}
\newcommand{\beqa}{\begin{eqnarray}}
\newcommand{\eeqa}[1]{\label{#1}\end{eqnarray}}
\newcommand{\bequ}{\begin{equation}}
\newcommand{\eequ}[1]{\label{#1}\end{equation}}
\newcommand{\Grad}{\nabla}
\newcommand{\Div}{\nabla \cdot}
\newcommand{\Curl}{\nabla \times}
\newcommand{\mn}{{\rm min}}
\newcommand{\mx}{{\rm max}}
\newcommand{\lang}{\langle}
\newcommand{\rang}{\rangle}
\newcommand{\Md}{\partial}
\newcommand{\MP}{\bigoplus}
\newcommand{\MT}{\bigotimes}
\newcommand{\MU}{\bigcup}
\newcommand{\MI}{\bigcap}
\newcommand{\MV}{\bigvee}
\newcommand{\MW}{\bigwedge}
\newcommand{\ov}[1]{\overline{#1}}

\newcommand{\Ga}{\alpha}
\newcommand{\Gb}{\beta}
\newcommand{\Gd}{\delta}
\newcommand{\Ge}{\epsilon}
\newcommand{\Gve}{\varepsilon}
\newcommand{\Gf}{\phi}
\newcommand{\Gvf}{\varphi}
\newcommand{\Gg}{\gamma}
\newcommand{\Gc}{\chi}
\newcommand{\Gi}{\iota}
\newcommand{\Gk}{\kappa}
\newcommand{\Gl}{\lambda}
\newcommand{\Gn}{\eta}
\newcommand{\Gm}{\mu}
\newcommand{\Gv}{\nu}
\newcommand{\Gp}{\pi}
\newcommand{\Gt}{\theta}
\newcommand{\Gvt}{\vartheta}
\newcommand{\Gr}{\rho}
\newcommand{\Gvr}{\varrho}
\newcommand{\Gs}{\sigma}
\newcommand{\Gvs}{\varsigma}
\newcommand{\Gj}{\tau}
\newcommand{\Gu}{\upsilon}
\newcommand{\Go}{\omega}
\newcommand{\Gx}{\xi}
\newcommand{\Gy}{\psi}
\newcommand{\Gz}{\zeta}
\newcommand{\GD}{\Delta}
\newcommand{\GF}{\Phi}
\newcommand{\GG}{\Gamma}
\newcommand{\GL}{\Lambda}
\newcommand{\GP}{\Pi}
\newcommand{\GT}{\Theta}
\newcommand{\GS}{\Sigma}
\newcommand{\GU}{\Upsilon}
\newcommand{\GO}{\Omega}
\newcommand{\GX}{\Xi}
\newcommand{\GY}{\Psi}

\newcommand{\BGa}{\bfm\alpha}
\newcommand{\BGb}{\bfm\beta}
\newcommand{\BGd}{\bfm\delta}
\newcommand{\BGe}{\bfm\epsilon}
\newcommand{\BGve}{\bfm\varepsilon}
\newcommand{\BGf}{\bfm\phi}
\newcommand{\BGvf}{\bfm\varphi}
\newcommand{\BGg}{\bfm\gamma}
\newcommand{\BGc}{\bfm\chi}
\newcommand{\BGi}{\bfm\iota}
\newcommand{\BGk}{\bfm\kappa}
\newcommand{\BGl}{\bfm\lambda}
\newcommand{\BGn}{\bfm\eta}
\newcommand{\BGm}{\bfm\mu}
\newcommand{\BGv}{\bfm\nu}
\newcommand{\BGp}{\bfm\pi}
\newcommand{\BGt}{\bfm\theta}
\newcommand{\BGvt}{\bfm\vartheta}
\newcommand{\BGr}{\bfm\rho}
\newcommand{\BGvr}{\bfm\varrho}
\newcommand{\BGs}{\bfm\sigma}
\newcommand{\BGvs}{\bfm\varsigma}
\newcommand{\BGj}{\bfm\tau}
\newcommand{\BGu}{\bfm\upsilon}
\newcommand{\BGo}{\bfm\omega}
\newcommand{\BGx}{\bfm\xi}
\newcommand{\BMx}{\bfm\xi}
\newcommand{\BGy}{\bfm\psi}
\newcommand{\BGz}{\bfm\zeta}
\newcommand{\BGD}{\bfm\Delta}
\newcommand{\BGF}{\bfm\Phi}
\newcommand{\BGG}{\bfm\Gamma}
\newcommand{\BGL}{\bfm\Lambda}
\newcommand{\BGP}{\bfm\Pi}
\newcommand{\BGT}{\bfm\Theta}
\newcommand{\BGS}{\bfm\Sigma}
\newcommand{\BGU}{\bfm\Upsilon}
\newcommand{\BGO}{\bfm\Omega}
\newcommand{\BGX}{\bfm\Xi}
\newcommand{\BGY}{\bfm\Psi}

\newcommand{\Ca}{{\cal a}}
\newcommand{\Cb}{{\cal b}}
\newcommand{\Cc}{{\cal c}}
\newcommand{\Cd}{{\cal d}}
\newcommand{\Ce}{{\cal e}}
\newcommand{\Cf}{{\cal f}}
\newcommand{\Cg}{{\cal g}}
\newcommand{\Ch}{{\cal h}}
\newcommand{\Ci}{{\cal i}}
\newcommand{\Cj}{{\cal j}}
\newcommand{\Ck}{{\cal k}}
\newcommand{\Cl}{{\cal l}}
\newcommand{\Cm}{{\cal m}}
\newcommand{\Cn}{{\cal n}}
\newcommand{\Co}{{\cal o}}
\newcommand{\Cp}{{\cal p}}
\newcommand{\Cq}{{\cal q}}
\newcommand{\Cr}{{\cal r}}
\newcommand{\Cs}{{\cal s}}
\newcommand{\Ct}{{\cal t}}
\newcommand{\Cu}{{\cal u}}
\newcommand{\Cv}{{\cal v}}
\newcommand{\Cx}{{\cal x}}
\newcommand{\Cy}{{\cal y}}
\newcommand{\Cz}{{\cal z}}
\newcommand{\CA}{{\cal A}}
\newcommand{\CB}{{\cal B}}
\newcommand{\CC}{{\cal C}}
\newcommand{\CD}{{\cal D}}
\newcommand{\CE}{{\cal E}}
\newcommand{\CF}{{\cal F}}
\newcommand{\CG}{{\cal G}}
\newcommand{\CH}{{\cal H}}
\newcommand{\CI}{{\cal I}}
\newcommand{\CJ}{{\cal J}}
\newcommand{\CK}{{\cal K}}
\newcommand{\CL}{{\cal L}}
\newcommand{\CM}{{\cal M}}
\newcommand{\CN}{{\cal N}}
\newcommand{\CO}{{\cal O}}
\newcommand{\CP}{{\cal P}}
\newcommand{\CQ}{{\cal Q}}
\newcommand{\CR}{{\cal R}}
\newcommand{\CS}{{\cal S}}
\newcommand{\CT}{{\cal T}}
\newcommand{\CU}{{\cal U}}
\newcommand{\CV}{{\cal V}}
\newcommand{\CW}{{\cal W}}
\newcommand{\CX}{{\cal X}}
\newcommand{\CY}{{\cal Y}}
\newcommand{\CZ}{{\cal Z}}
\newcommand{\BCa}{{\bfm{\cal a}}}
\newcommand{\BCb}{{\bfm{\cal b}}}
\newcommand{\BCc}{{\bfm{\cal c}}}
\newcommand{\BCd}{{\bfm{\cal d}}}
\newcommand{\BCe}{{\bfm{\cal e}}}
\newcommand{\BCf}{{\bfm{\cal f}}}
\newcommand{\BCg}{{\bfm{\cal g}}}
\newcommand{\BCh}{{\bfm{\cal h}}}
\newcommand{\BCi}{{\bfm{\cal i}}}
\newcommand{\BCj}{{\bfm{\cal j}}}
\newcommand{\BCk}{{\bfm{\cal k}}}
\newcommand{\BCl}{{\bfm{\cal l}}}
\newcommand{\BCm}{{\bfm{\cal m}}}
\newcommand{\BCn}{{\bfm{\cal n}}}
\newcommand{\BCo}{{\bfm{\cal o}}}
\newcommand{\BCp}{{\bfm{\cal p}}}
\newcommand{\BCq}{{\bfm{\cal q}}}
\newcommand{\BCr}{{\bfm{\cal r}}}
\newcommand{\BCs}{{\bfm{\cal s}}}
\newcommand{\BCt}{{\bfm{\cal t}}}
\newcommand{\BCu}{{\bfm{\cal u}}}
\newcommand{\BCv}{{\bfm{\cal v}}}
\newcommand{\BCx}{{\bfm{\cal x}}}
\newcommand{\BCy}{{\bfm{\cal y}}}
\newcommand{\BCz}{{\bfm{\cal z}}}
\newcommand{\BCA}{{\bfm{\cal A}}}
\newcommand{\BCB}{{\bfm{\cal B}}}
\newcommand{\BCC}{{\bfm{\cal C}}}
\newcommand{\BCD}{{\bfm{\cal D}}}
\newcommand{\BCE}{{\bfm{\cal E}}}
\newcommand{\BCF}{{\bfm{\cal F}}}
\newcommand{\BCG}{{\bfm{\cal G}}}
\newcommand{\BCH}{{\bfm{\cal H}}}
\newcommand{\BCI}{{\bfm{\cal I}}}
\newcommand{\BCJ}{{\bfm{\cal J}}}
\newcommand{\BCK}{{\bfm{\cal K}}}
\newcommand{\BCL}{{\bfm{\cal L}}}
\newcommand{\BCM}{{\bfm{\cal M}}}
\newcommand{\BCN}{{\bfm{\cal N}}}
\newcommand{\BCO}{{\bfm{\cal O}}}
\newcommand{\BCP}{{\bfm{\cal P}}}
\newcommand{\BCQ}{{\bfm{\cal Q}}}
\newcommand{\BCR}{{\bfm{\cal R}}}
\newcommand{\BCS}{{\bfm{\cal S}}}
\newcommand{\BCT}{{\bfm{\cal T}}}
\newcommand{\BCU}{{\bfm{\cal U}}}
\newcommand{\BCV}{{\bfm{\cal V}}}
\newcommand{\BCW}{{\bfm{\cal W}}}
\newcommand{\BCX}{{\bfm{\cal X}}}
\newcommand{\BCY}{{\bfm{\cal Y}}}
\newcommand{\BCZ}{{\bfm{\cal Z}}}

\def\ii{{\rm i}}
\def\dd{{\rm d}}
\def\Im{{\it Im}}
\def\Re{{\it Re}}
\def\ca{{\cal A}}
\def\ct{{\cal T}}
\def\Ba{{\bf a}}
\def\Bb{{\bf b}}
\def\Bc{{\bf c}}
\def\Bd{{\bf d}}
\def\Be{{\bf e}}
\def\Bf{{\bf f}}
\def\Bg{{\bf g}}
\def\Bh{{\bf h}}
\def\Bi{{\bf i}}
\def\Bj{{\bf j}}
\def\Bk{{\bf k}}
\def\Bl{{\bf l}}
\def\Bm{{\bf m}}
\def\Bn{{\bf n}}
\def\Bo{{\bf o}}
\def\Bp{{\bf p}}
\def\Bq{{\bf q}}
\def\Br{{\bf r}}
\def\Bs{{\bf s}}
\def\Bt{{\bf t}}
\def\Bu{{\bf u}}
\def\Bv{{\bf v}}
\def\Bw{{\bf w}}
\def\Bx{{\bf x}}
\def\By{{\bf y}}
\def\Bz{{\bf z}}
\def\BA{{\bf A}}
\def\BB{{\bf B}}
\def\BC{{\bf C}}
\def\BD{{\bf D}}
\def\BE{{\bf E}}
\def\BF{{\bf F}}
\def\BG{{\bf G}}
\def\BH{{\bf H}}
\def\BI{{\bf I}}
\def\BJ{{\bf J}}
\def\BK{{\bf K}}
\def\BL{{\bf L}}
\def\BM{{\bf M}}
\def\BN{{\bf N}}
\def\BO{{\bf O}}
\def\BP{{\bf P}}
\def\BQ{{\bf Q}}
\def\BR{{\bf R}}
\def\BS{{\bf S}}
\def\BT{{\bf T}}
\def\BU{{\bf U}}
\def\BV{{\bf V}}
\def\BW{{\bf W}}
\def\BX{{\bf X}}
\def\BY{{\bf Y}}
\def\BZ{{\bf Z}}

\def\half{{\scriptstyle{1\over 2}}}
\newcommand{\eps}{\varepsilon}
\newcommand{\beq}{\begin{equation}}
\newcommand{\eeq}{\end{equation}}
\newcommand{\overliner}{\begin{eqnarray}}
\newcommand{\earr}{\end{eqnarray}}
\newcommand{\beqn}{\begin{equation*}}
\newcommand{\eeqn}{\end{equation*}}
\newcommand{\overlinern}{\begin{eqnarray*}}
\newcommand{\earrn}{\end{eqnarray*}}
\newcommand{\prt}{\partial}
\newcommand{\sg}{\sigma}
\newcommand{\fr}{\frac}
\newcommand{\kap}{\varkappa}
\newcommand{\diag}{\mbox{diag}}
\newcommand{\sign}{\mbox{sign}}
\def\l{\label}

\title{\huge Corrugation crack front waves}

\author{J.R. Willis$^1$, N.V. Movchan$^2$ and A.B. Movchan$^2$ \\
$^1$ {\it Department of Applied Mathematics
        and Theoretical Physics, }\\ {\it University of Cambridge,
       Cambridge, U.K.} \\
       $^2$ {\it Department of Mathematical Sciences, University of Liverpool, }
\\
       {\it     Liverpool, U.K.}
}

\date{}
\maketitle

\begin{abstract}
The paper presents a model of a dynamic crack with a wavy surface.
So far, theoretical analysis of crack front waves has been performed only
for in-plane perturbations of the crack front. In the present paper,
generalisation is given to a more general three-dimensional perturbation,
and equations that govern corrugation crack front waves are derived
and analysed.

\end{abstract}

\noindent {\bf Keywords:} Dynamic fracture, crack front waves, asymptotic
analysis.

\section{Introduction}

The paper analyses  singular fields around a dynamic
crack whose surface is slightly perturbed from the original
plane configuration.  Crack front waves in the plane of the crack
were discovered numerically by Morrissey and Rice in \cite{MR}, and
later  confirmed analytically  by Ramanathan and
Fisher \cite{RF}, using the results of the perturbation
analysis of Willis and Movchan \cite{WM1}.
Experimental observations of persistent crack front
waves were reported by Sharon, Cohen and Fineberg \cite{SCF}.
The more general development of Willis and Movchan \cite{WM_visco} and
Woolfries {\em et al.} \cite{WMW}
extended the analysis to a crack propagating through a viscoelastic
medium. The perturbation formulae for the stress intensity factors,
specialised to a plane strain formulation, have been used by Obrezanova
{\em et al.} \cite{OMW1} in the stability analysis of rectilinear
propagation. A quasi-static advance of a tunnel crack under a mixed mode loading has been analysed by 
Lazarus and Leblond \cite{LL}. 

The aim of the present paper is to develop a model describing
corrugation (out-of-plane) waves along the front of a moving crack.
This work is based on the ideas of the earlier publication by
Willis \cite{W_corr}. The plan of the paper is as follows.
We begin, in Section 2,
with the description of the geometry, governing equations and
perturbation functions. A summary of the first-order
approximations for the stress intensity factors is presented in Section 2.2.
Section 3 includes the study of the corrugation waves in the
first-order asymptotic approximation for a basic Mode I loading.
In Section 4, we derive the dispersion equation for crack front waves in the mixed mode I-III loading.
The technical appendix contains an outline of the fundamental integral
identity, and the expressions for effective tractions. 


\section{Basic perturbation formulae}

For a linearly elastic medium, we consider a semi-infinite crack with
a slightly perturbed surface.
The unperturbed configuration of the crack at time $t$ is defined by
\begin{equation}
S_0(t) = \{{\bf x}: \; -\infty < x_1 < Vt,
\; -\infty < x_2 < \infty, \; x_3 =0\},
\end{equation}
where $V$ is a constant crack speed, which does not exceed the Rayleigh
wave speed.
The perturbation is introduced through deviations of the crack front
in both in-plane and out-of-plane directions. The
perturbed surface of the crack at time $t$ is
\begin{eqnarray}
S_\varepsilon(t) &=& \{{\bf x}: - \infty < x_1 < Vt +\varepsilon\Gvf(x_2,t),
\nonumber\\
&& -\infty < x_2 < \infty,\; x_3 = \varepsilon
\psi(x_1-Vt,x_2)  
 \}. 
\label{Seps}
\end{eqnarray}
The functions $\Gvf$ 
and $\psi$ 
are smooth and bounded, and $\Gve$ is a small
non-dimensional parameter,
$0\leq \varepsilon \ll 1$.
It is helpful to use the moving-frame coordinates, so that
$X=x_1-Vt$.

It is assumed that the medium is loaded
so that a stress $\BGs^{\rm nc}$ and a displacement
${\bf u}^{\rm nc}$ would be generated in the absence of the crack.
The crack induces additional fields $\BGs$, ${\bf u}$.
They
satisfy the equations of motion and the traction boundary
conditions on the crack faces:
\beq
\Gs_{ij,j} - \rho \ddot u_i=0,~ i=1,2,3, ~~ \hbox{outside the crack}
\eequ{eqmotion}
and
\begin{equation}
\sigma_{ij}n_j + \sigma_{ij}^{\rm nc}n_j = 0, ~ \hbox{ on the crack faces},
\label{3.2w}\end{equation}
and 
correspond to waves outgoing from the crack as $x_3 \rightarrow \pm \infty.$

\subsection{Local coordinates and asymptotics for stresses}

At a point $\Bx^0=(x_1^0,x_2^0,x_3^0),$ which is on the crack edge at time $t,$ so that
$$
x_1^0=Vt + \varepsilon \varphi(x_2^0,t), 
\; x_3^0=\varepsilon \psi(x_1^0-Vt, x_2^0), 
$$
we define a coordinate system such that
$$
\Bx-\Bx^0= \sum_{i=1}^3 x'_i \Be'_i,
$$
where
\beq
\pmatrix{\Be'_1 \cr \Be'_2 \cr \Be'_3}
= \Bigg\{ \BI + \Gve\left(\begin{array}{ccc}
        0 & -\Gvf_{,2} & \psi^*_{,1}\\
        \Gvf_{,2} & 0 &  \psi^*_{,2}\\
         - \psi^*_{,1} & - \psi^*_{,2} & 0
  \end{array}\right)
%
\Bigg\}
\pmatrix{\Be_1 \cr \Be_2 \cr \Be_3}.
\label{3.3w}\end{equation}
Here $\psi^*$ denotes $\psi$ evaluated for $x_1=Vt.$
The above transformation involves a shift to the crack edge and a further
rotation of coordinate axes.

\noindent
In the new frame, the stress components ($\sigma_{i3}^\prime$) have the asymptotic
form
$$
\sigma_{i3}^\prime(x'_1, x'_2,0)\sim (K_i^{(0)} + \varepsilon K_i^{(1)} 
)/(2\pi
x_1^\prime)^{1/2}
 - (P_i^{(0)} + \GD P_i)
$$
\begin{equation}
+
(A_i^{(0)} + \GD A_i)(x_1^\prime)^{1/2} - (F_i^{(0)} + \GD F_i) x_1^\prime
+(N_i^{(0)} + \GD N_i) (x_1^\prime)^{3/2}, ~i=1,2,3.
\label{3.4w}\end{equation}

The first-order asymptotic approximation of stress-intensity factors
was constructed and studied
in \cite{WM1}, \cite{MW2}, \cite{WM97}, \cite{W}.
In Appendix  we include a description of the fundamental
identity, which is essential for this work. We also require the dynamic
crack face weight function $[\BU]$, as defined in Appendix. The field $[{\bf U}]$ has a singularity proportional to $X^{-1/2}H(X)\delta(x_2)\delta(t)$ as
$X\to 0$.

\subsection{First-order perturbations of the stress intensity factors}

We begin with the first-order approximation for the stress intensity factors,
when
\begin{equation}
{K}_j \sim {K}_j^{(0)} + \varepsilon{K}_j^{(1)}, ~ j=I, II, III. \label{SIF_exp}
\end{equation}
For the Mode-I unperturbed case, $K_{II}^{(0)}= K_{III}^{(0)}=0$, and the
perturbation terms are defined by (see \cite{WM97}, \cite{W})
\begin{eqnarray}
K_{II}^{(1)}&=& - Q_{11}*\psi^*\Theta_{13}K_I^{(0)}
- \psi^*_{,1}\omega_{13}K_I^{(0)}
- \psi^* \left(
\Sigma_{11} +\frac{V^2}{2b^2}\Sigma_{12}\right) A_3^{(0)} \sqrt{\fr{\pi}{2}}
\nonumber\\
&& \hbox{\hskip .3in}+ [U]_{11}*\langle P_1^{(1)}\rangle
+ [U]_{21}*\langle P_2^
{(1)}\rangle
 -\langle U\rangle_{31}*[P_3^{(1)}],
\label{3.7w}\\
K_{III}^{(1)} &=& -Q_{12}*\psi^*
\Theta_{13} K_I^{(0)} - \psi^*_{,2}\omega_{23}K_I^{(0)}
\nonumber\\
&& \hbox{\hskip .3in}+[U]_{12}*\langle P_1^{(1)}\rangle + [U]_{22}*\langle P_2^{
(1)}\rangle
 -\langle U\rangle_{32}*[P_3^{(1)}],
\label{3.8w}\\
K_I^{(1)} &=& Q_{33}*\Gvf K_I^{(0)}
+\left(\frac{\pi}{2}\right)^{1/2} \Gvf A_3^{(0)}
 - \langle U \rangle_{13}*[P_1^{(1)}]
\nonumber\\
&&
\hbox{\hskip .3in}
-
\langle U \rangle_{23}*[ P_2^
{(1)}]+[U]_{33}*\langle P_3^{(1)}\rangle.
\label{3.9w}\end{eqnarray}
The matrix $\BQ$ is a block-diagonal matrix defined in \cite{WM97}; other
functions that appear in the above equations are
$$
\Theta_{13} = \Sigma_{11}+\frac{V^2}{2b^2}\Sigma_{12},
~ \omega_{13} = \frac{\alpha-\beta}{R(V)}(1+\beta^2)(\alpha +2\beta) -2,
$$
$$
\omega_{23} = \frac{2\nu}{R(V)}(1+\beta^2)(\alpha^2-\beta^2) -1,
$$
$$
\Sigma_{11} = -\frac{4\alpha\beta -(1+2\alpha^2-\beta^2)(1+\beta^2)}{R(V)}, ~
\Sigma_{12} = \frac{-2(1+\beta^2 - 2\alpha\beta)}{R(V)},
$$
\beq
\alpha^2 = 1-V^2/a^2,\;\;\;\beta^2 = 1-V^2/b^2, ~
R(V) = 4\alpha\beta-(1+\beta^2)^2.
\eequ{3.10w}
Here, $a$ and $b$ denote the speeds of longitudinal and shear
waves, respectively. The representations for the effective tractions $P_i^{(1)},
i=1,2,3,$ are given in Appendix.

\subsection{Crack front waves confined to the plane $x_3=0$}
\label{inplane}

Assuming that the out-of-plane deflection is not present ($\psi=0$),
we consider a  first-order  
in-plane perturbation of the crack front and loading in Mode I, so that $\sigma_{13}^{\rm nc} =  \sigma_{23}^{\rm nc} = 0$ on the plane $x_3=0$. 
In this special case,
the only non-zero
stress intensity factor is $K_I$, and the corresponding
perturbation formula
reduces to
\begin{equation}
K_I^{(1)} = Q_{33}*\Gvf K_I^{(0)} +\left(\frac{\pi}{2}\right)^{1/2}\Gvf
A_3^{(0)}.
\label{4.1}\end{equation}

According to the Griffith energy balance equation,
the energy flux ${\cal G}$
into the crack edge is constant, denoted here by ${\cal G}_c$:
\begin{equation}
{\cal G}\equiv \frac{1- \nu^2}{E} f_I(v) K_I^2 = {\cal G}_c.
\label{4.2}\end{equation}
Here, $v$ is the local crack speed (to the first-order approximation,
$v=V+ \Gve \dot \Gvf$)
and $f_I(v)$ is a known function (e.g., \cite{F}):
\begin{equation}
f_{I}(v)
=\frac{v^2 \; \alpha(v)}{(1-\nu) \; b^2 \; R(v)}.
\label{4.3}\end{equation}
Expanding the Griffith energy balance equation
(\ref{4.2}) to order $\varepsilon$,
we obtain
\begin{equation}
2Q_{33}*\Gvf +\frac{f^\prime_I(V)}{f_I(V)}\dot\Gvf + 2 m\Gvf = 0,
\label{4.4}\end{equation}
where $m=(\pi/2)^{1/2}A_3^{(0)}/K_I^{(0)}$.
Applying the Fourier transform with respect to $t$ and $x_2$ we deduce
that
a non-zero solution is possible only if the dispersion relation
\begin{equation}
2\overline Q_{33}(\omega,k) -{\rm i}\omega\frac{f^\prime_I(V)}{f_I(V)} + 2m = 0
\label{4.5}\end{equation}
is satisfied. Here, the Fourier transform
$\overline Q_{33}$ is a homogeneous function of
degree 1 in $(\omega, k)$. At high
frequency and large wavenumber, the third term in the above equation
can be neglected. Such an equation can be solved for $\omega/k$, and a real
root represents a speed of wave propagating along the crack front.
This computation was performed by Ramanathan and Fisher \cite{RF}.

\section{Corrugation waves for a Mode-I basic loading. First-order analysis.}

Can a Mode-I basic loading generate a corrugation wave propagating along the
crack front? This case corresponds to a non-zero out-of-plane perturbation
characterised by the function $\psi(x_1-Vt, x_2)$. 
Crack stability with respect to out-of-plane deflections
can be studied, once a fracture criterion is
identified.

\label{faa}

If we suppose that $K_{II}=0$ then, to lowest order, $\psi$
must satisfy $K_{II}^{(1)}(\psi) =0$,
where $K_{II}^{(1)}$ is given by \eq{3.7w}.
The proposition that the crack propagates so as to
maintain $K_{II} = 0$ together with the Griffith
energy balance has recently received theoretical support,
on the basis of a version of Hamilton's principle \cite{Oleaga}.

\begin{figure}[h]
\centerline{\psfig{file=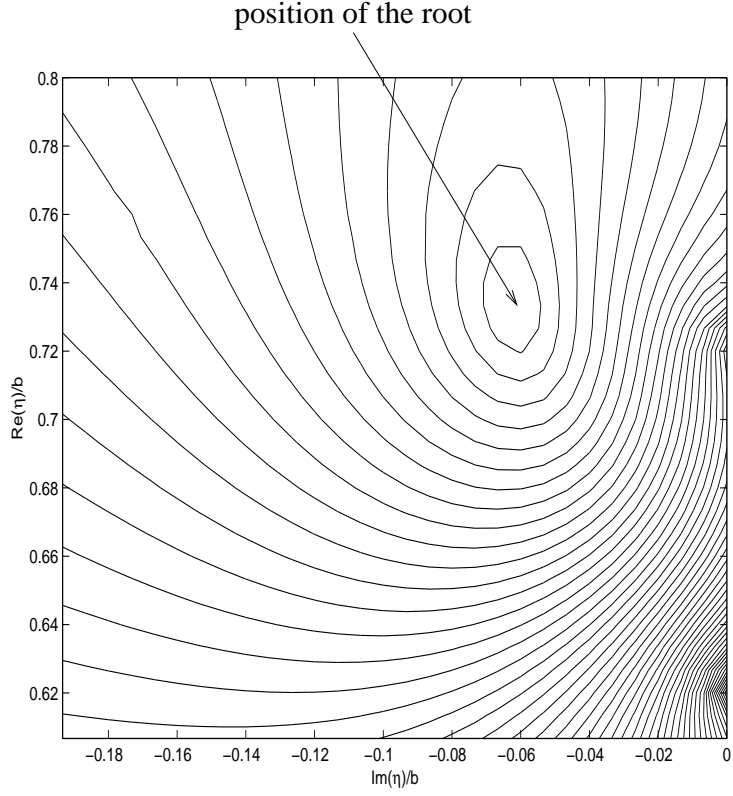,height=10.5cm}}
\caption{Level curves of the function ${\cal W}=|\overline Q_{11}\Theta_{13}  - {\rm i} (\omega/V) \omega_{13}|$, for
$V/b=0.69$ and $\nu=0.3$.
}
\label{diagram}
\end{figure}

Assuming that the in-plane perturbation of the crack front equals zero,
we look into stability against out-of-plane deflections.
It is also assumed that $\omega\equiv k_1 V$ and $k_2$
are large.
The leading-order approximation of the stress intensity factor
$K_{II}$ yields
\begin{equation}
\overline K_{II}^{(1)}=  \{-\overline Q_{11}\Theta_{13}  +
{\rm i} (\omega/V) \omega_{13}\}K_I^{(0)}\overline\psi^*=0. \label{5.1}
\end{equation}
This relation is homogeneous of degree 1 in $\omega$ and $k_2$, and so is
non-dispersive.

The numerical study of equation \eq{5.1} produced the following results.


$\bullet$ For crack speeds $V$ greater than a critical value $V_c$ (which is close to $0.6$
of the Rayleigh wave speed) there is a value $\eta = \omega / |k_2|$ with small, negative,
imaginary part that satisfies (\ref{5.1}). The position of the root is shown in Figure
\ref{diagram}; the calculation is produced for the case of $V/b=0.69$, and the diagram shows
the level curves of the modulus of the expression in the curly brackets on
the left side of (\ref{5.1}).

\begin{figure}[h]
\centerline{\psfig{file=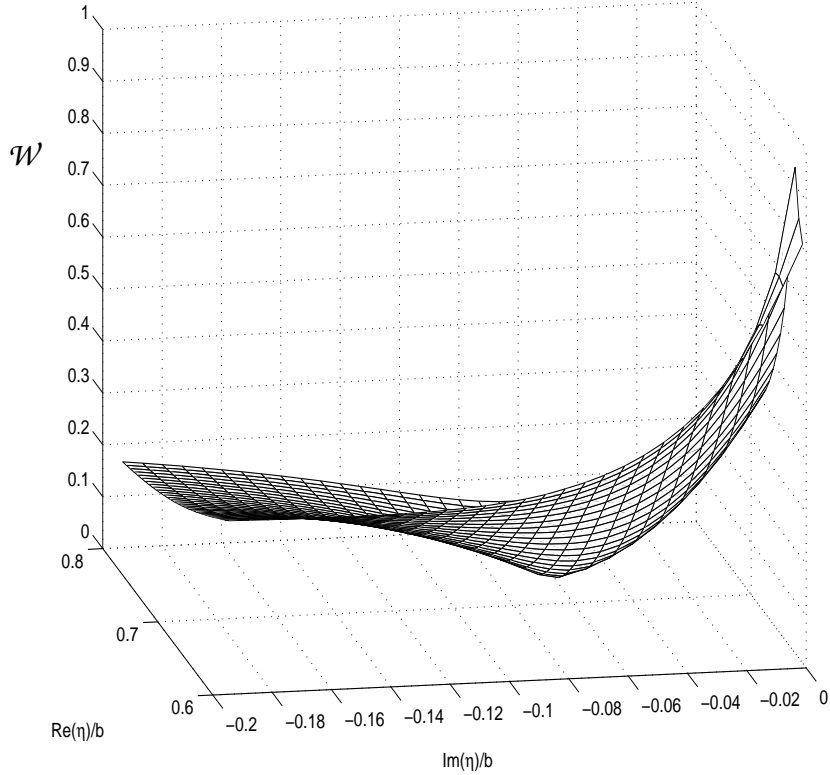,height=10.5cm}
}
\caption{Surface plot of the function ${\cal W}=|\overline Q_{11} \Theta_{13} - {\rm i}
(\omega/V) \omega_{13}|$, for $V/b=0.69$ and $\nu=0.3.$
}
\label{diagram1}
\end{figure}

Figure \ref{diagram} is accompanied by a three dimensional surface plot, shown in Figure \ref{diagram1},  of the 
function 
${\cal W}=|\overline Q_{11}\Theta_{13}  - {\rm i} (\omega/V) \omega_{13}|$; the surface touches
the $\eta$-plane at the point corresponding to the root of equation (\ref{5.1}).

\begin{figure}[h]
\centerline{\psfig{file=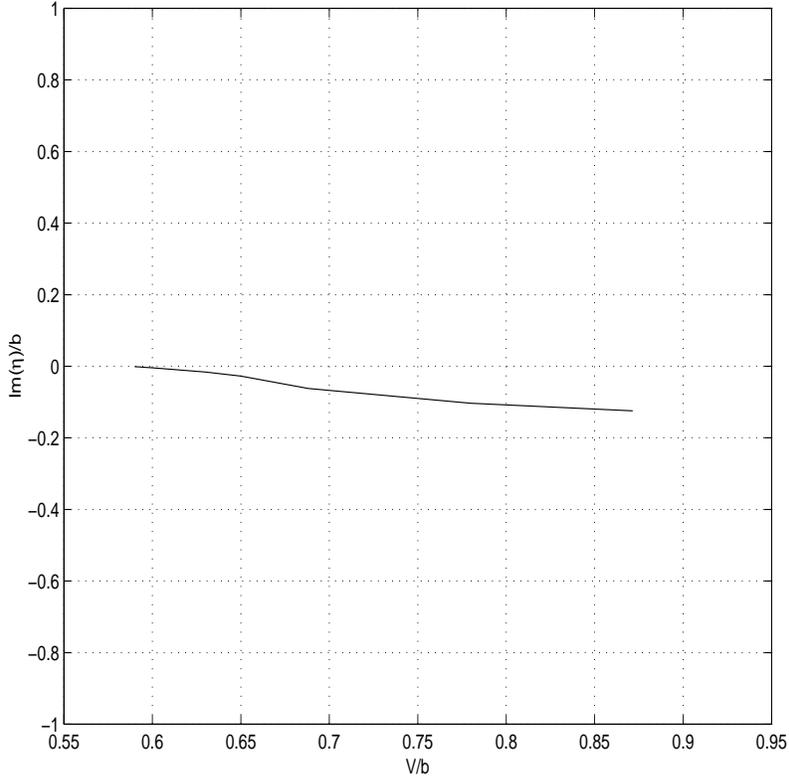,height=10.5cm}
}
\caption{The imaginary part of $\eta$ as a function of the crack front velocity $V$.
}
\label{diagram2}
\end{figure}

$\bullet$ The "corrugation wave" suffers slow attenuation as it propagates. The imaginary part
of $\eta$, which characterises the rate of attenuation of the
"corrugation wave", is shown in
Fig. \ref{diagram2} for different values of the crack front velocity $V$, and it decreases
with $V$.

\section{First-order coupling between in-plane and out-of-plane crack
front perturbations for mixed Mode I-III loading}

Here, we assume that $K_{II}^{(0)}=0$, whereas $K_{I}^{(0)}$  and $K_{III}^{(0)}$
are non-zero for a half-plane crack propagating with constant
speed $V$ (unperturbed configuration).
To  first order, the stress intensity factors
are represented by the formulae (\ref{SIF_exp}),
where the
perturbation terms ${K}_j^{(1)}, ~ j=I, II, III,$ are defined by (see \cite{WM97}, \cite{W})
\begin{eqnarray}
K_{II}^{(1)}&=& - Q_{11}*\psi^*\Theta_{13}K_I^{(0)}
- \psi^*_{,1}\omega_{13}K_I^{(0)}
- \psi^* \left(
\Sigma_{11} +\frac{V^2}{2b^2}\Sigma_{12}\right) A_3^{(0)} \sqrt{\fr{\pi}{2}}
\nonumber\\
&& +Q_{21} * (\Gvf K_{III}^{(0)}) - \Gvf_{,2} K_{III}^{(0)} + \sqrt{\fr{\pi}{2}} \Gvf
A_1^{(0)} \nonumber \\
&&
+ [U]_{11}*\langle P_1^{(1)}\rangle
+ [U]_{21}*\langle P_2^
{(1)}\rangle
 -\langle U\rangle_{31}*[P_3^{(1)}],
\label{3.7wa}\\
K_{III}^{(1)} &=& -Q_{12}*\psi^*
\Theta_{13} K_I^{(0)} - \psi^*_{,2}\omega_{23}K_I^{(0)}
+ Q_{22}* (\Gvf K_{III}^{(0)}) + \sqrt{\fr{\pi}{2}} \Gvf A_{2}^{(0)}
\nonumber\\
&&
+[U]_{12}*\langle P_1^{(1)}\rangle + [U]_{22}*\langle P_2^{
(1)}\rangle
 -\langle U\rangle_{32}*[P_3^{(1)}],
\label{3.8wa}\\
K_I^{(1)} &=& Q_{33}*\Gvf K_I^{(0)}
+\left(\frac{\pi}{2}\right)^{1/2} \Gvf A_3^{(0)}
- \psi^*\Big(1 - \fr{V^2}{2 b^2} \GS_{12} \Big) A_1^{(0)}   \sqrt{\fr{\pi}{2}}
\nonumber\\
&&
- 2 \psi^*_{,2} K_{III}^{(0)} - \langle U \rangle_{13}*[P_1^{(1)}] -
\langle U \rangle_{23}*[ P_2^
{(1)}]+[U]_{33}*\langle P_3^{(1)}\rangle.
\label{3.9wa}\end{eqnarray}

We shall use the criterion of local symmetry
$
K_{II}=0,
$
together with the
Griffith energy balance equation
\begin{equation}
{\cal G}\equiv(2 \mu)^{-1} f_I(v) K_I^2 + (2 \mu)^{-1} f_{III}(v) K_{III}^2
={\cal G}_c=\mbox{const}. \label{EBeqn}
\end{equation}

Taking into account that, to first order,
$
v \sim V+ \varepsilon \dot{\varphi}
$, we deduce
$$
{\cal G}=(2 \mu)^{-1} f_I(V) (K_I^{(0)})^2 + (2 \mu)^{-1} f_{III}(V) (K_{III}^{(0)})^2
$$
$$
+ \varepsilon (2 \mu)^{-1} \Big(\dot{\varphi} f'_I(V)(K_I^{(0)})^2+2 f_I(V) K_I^{(0)} K_I^{(1)}
$$
\beq
+\dot{\varphi} f'_{III}(V)(K_{III}^{(0)})^2+2 f_{III}(V) K_{III}^{(0)} K_{III}^{(1)}\Big)+O(\varepsilon^2).
\eequ{EB1}
It follows from (\ref{EBeqn}), \eq{EB1}
and the local symmetry criterion
$K_{II}=0$
that
$$
\dot{\varphi} (f'_I(V)(K_I^{(0)})^2+f'_{III}(V)(K_{III}^{(0)})^2)
+2 f_I(V) K_I^{(0)} K_I^{(1)}(\varphi, \psi)
$$
\beq
+2f_{III}(V) K_{III}^{(0)} K_{III}^{(1)}(\varphi,\psi)=0, \ \ \
\eequ{energy_release}
\beq
K_{II}^{(1)}(\varphi, \psi)=0.
\eequ{KII1}
The above equations define the coupling between the in-plane and
out-of-plane perturbations of the crack front.

Applying the Fourier transform with respect to $t$ and $x_2$ and assuming that
$\omega=k_1 V$ and $k_2$ are large, we deduce
$$
\ov{\varphi}\Big\{\left(2f_I(V)\ov{Q}_{33}-\mbox{i} \omega f'_I(V)\right)(K_I^{(0)})^2+
\left(2f_{III}(V)\ov{Q}_{22}- \mbox{i} \omega f'_{III}(V)\right)(K_{III}^{(0)})^2  \Big\}
$$
\begin{equation}
+\ov{\psi^*}K_I^{(0)} K_{III}^{(0)} \Big\{2 f_{III}(V)
\left(-\ov{Q}_{12} \Theta_{13}+\mbox{i} k_2 \omega_{23}\right)+4 \mbox{i} k_2 f_I(V)
\Big\}=0, \label{FTeqn1}
\end{equation}
\begin{equation}
\{-\ov{Q}_{11} \Theta_{13}+\mbox{i} (\omega/V) \omega_{13}\}
\ov{\psi^*} K_I^{(0)} + (\ov{Q}_{21} + \mbox{i} k_2) \ov{\varphi} K_{III}^{(0)}=0. \label{FTeqn2}
\end{equation}
The system (\ref{FTeqn1}), (\ref{FTeqn2}) is linear in $\ov{\varphi}$ and $\ov{\psi^*}$,
and it possesses a nontrivial solution
if and only if
the matrix of this system is degenerate. This yields
the following dispersion relation:
$$
\Big\{-\ov{Q}_{11} \Theta_{13}+\mbox{i}(\omega/V) \omega_{13}\Big\}
\Big\{2f_I(V)\ov{Q}_{33}-\mbox{i} \omega f'_I(V)+
\left(2f_{III}(V)\ov{Q}_{22}- \mbox{i} \omega f'_{III}(V)\right){\cal K}_0^2  \Big\}
$$
\begin{equation}
-{\cal K}_0^2 \left(\ov{Q}_{21}+\mbox{i} k_2 \right)\Big\{2 f_{III}(V)
\left(-\ov{Q}_{12} \Theta_{13}+\mbox{i} k_2 \omega_{23}\right)+4 \mbox{i} k_2 f_I(V)
\Big\}=0. \label{dispeqn}
\end{equation}
Here ${\cal K}_0=K_{III}^{(0)}/K_I^{(0)}.$ The above dispersion equation, connecting $\omega$ and $k_2$, is to be analysed numerically to identify possible crack front waves associated with the 
external mixed mode I-III load. 

\section*{Appendix. Fundamental identity and effective tractions.}

\renewcommand{\theequation}{A\arabic{equation}}
\setcounter{equation}{0}

Here, we briefly describe the method developed in \cite{WM1},
\cite{MW2}, \cite{WM97}.
We use the relation
\begin{equation}
{\bf u} = -{\bf G}*{\BGs},
\label{2.1w}\end{equation}
where ${\bf u}$ and $\BGs$ denote the values of the displacement
vector $(u_i)$ and the traction vector
$(\sigma_{i3})$ on the surface $x_3=0$ of the half-space $x_3 >0$; ${\bf G}$
is the Green's matrix function.
The symbol
$*$ denotes convolution over $x_1$, $x_2$ and $t$. It is assumed that all
waves emanate from the surface $x_3=0$. A similar identity applies to the
half-space
$x_3 <0$, with ${\bf G}$ being replaced by
$-{\bf G}^T$.

Three column vectors like ${\bf u}$ can be written  side by
side
to form a matrix ${\bf U}(+0)$, and similarly ${\BGS}(+0)$
represents the matrix formed from the
three corresponding vectors ${\BGs}$. Then
\begin{equation}
{\bf U}(+0) = -{\bf G}*{\BGS}(+0).
\label{2.2w}\end{equation}
The argument $(+0)$ signifies  values on
the boundary of the upper half-space.
Applying similar reasoning to the identity for the lower half-space
$x_3 <0$ gives
\begin{equation}
{\bf U}(-0) = {\bf G}^T*{\BGS}(-0).
\label{2.3w}\end{equation}
Next, we note that
\begin{eqnarray}
\{{\bf U}(+0)\}^T *{\BGs}(-0) &=& - \{{\BGS}(+0)\}^T*{\bf G}^T*{\BGs}(-0)
\nonumber\\
      &=& -\{{\BGS}(+0)\}^T*{\bf u}(-0),\label{2.4w}\\
\{{\bf U}(-0)\}^T *{\BGs}(+0) &=&  \{{\BGS}(-0)\}^T*{\bf G}^T*{\BGs}(+0)
\nonumber\\
      &=& -\{{\BGS}(-0)\}^T*{\bf u}(+0).
\label{2.5w}\end{eqnarray}
Subtracting the second line from the first and rearranging gives the identity
\begin{equation}
[{\bf U}]^T*\langle{\BGs}\rangle - \langle{\bf U}\rangle^T * [{\BGs}]
= -[{\BGS}]^T*\langle{\bf u}\rangle +\langle{\BGS}\rangle^T*[{\bf u}],
\label{2.6w}\end{equation}
where $\langle f\rangle = \frac{1}{2}(f(+0) + f(-0))$ and $[f] = f(+0) - f(-0)$.

In the moving frame associated with the crack edge, we use the coordinate
$X = x_1 -Vt$.
The operation of convolution
survives, with functions regarded as functions of
$X, x_2, t$ and the convolutions taken over these new variables.

For the unperturbed crack problem,
\begin{equation}
[{\BGs}] \equiv 0,\;\; [{\bf u}] = 0 \hbox{ when } X>0,\;\;
{\BGs}\equiv \langle{\BGs}\rangle = -{\BGs}^{\rm nc}\hbox{ when } X<0.
\label{2.7w}\end{equation}
We interpret equation (\ref{2.6w}) relative to the moving frame, and perform
factorizations of the Green's function so that ${\bf U}$ and ${\BGS}$ display
the related properties
\begin{equation}
[{\BGS}]\equiv 0,\;\; [{\bf U}] = 0 \hbox{ when } X<0,\;\;
{\BGS} \equiv \langle{\BGS}\rangle = 0 \hbox{ when } X>0.
\label{2.8w}\end{equation}
Equations \eq{2.2w}, \eq{2.3w} yield
\begin{equation}
[{\bf U}] = -({\bf G} + {\bf G}^T)*\langle{\BGS}\rangle,\;\;
\langle{\bf U}\rangle = -{\scriptstyle\frac{1}{2}}({\bf G}
-{\bf G}^T)*\langle{\BGS}\rangle.
\label{2.9w}\end{equation}
The first of these relations defines a Wiener--Hopf problem;
the second then gives $\langle{\bf U}\rangle$ directly. The Wiener--Hopf
problem uncouples into two sub-problems. One, associated with
the opening mode I of the crack, is a scalar problem. It was solved
in the case of elasticity in \cite{WM1}, and for a viscoelastic medium
in \cite{WMW}.
The remaining problem involves modes II and III, coupled. It was solved in
\cite{MW2}.

The field $[{\bf U}]$
has a singularity proportional to $X^{-1/2}H(X)\delta(x_2)\delta(t)$ as
$X\to 0$. With the constant of proportionality
chosen as $(2/\pi)^{1/2}{\bf I}$, we
call $[{\bf U}]$ the {\it dynamic weight function} for
the crack problem. With this choice,
letting $X\to +0$ in the identity (\ref{2.6w})
generates
\begin{equation}
{\bf K} = \lim_{X\to +0}\Big\{
\langle \BU \rangle^T * [\BGs^{(0)}]
-[{\bf U}]^T* \langle \BGs^{(0)} \rangle
\Big\},
\label{2.10w}\end{equation}
where ${\bf K}$ denotes the vector of stress-intensity factors $(K_{II},
K_{III}, K_I)^T$.
The matrix
function $\langle{\bf U}\rangle$ represents
a dynamical version of Bueckner's non-symmetric weight function,
as described in \cite{B} and \cite{MW2}.

We assume that the unperturbed steady-state crack is subjected to a
Mode-I loading, and the unperturbed displacement field is a vector function
$\Bu^{(0)} = \Bu^{(0)}(x_1-V t, x_2, x_3)$. We can write the resulting
displacement field in the form
$$
\Bu \sim \Bu^{(0)} + \Gve \Bu^{(1)}, 
$$
where $\Gve$ is a perturbation parameter.

The effective tractions
$P_{i}^{(1)}:=$
$-\Gs_{i3}(\Bu^{(1)})|_{x_3=0}, ~ i=1,2,3,$
have the form
(see formula (4.11) of \cite{WM97})
$$
P_i^{(1)} = -
\sum_{k=1}^2(\psi \Gs_{ik}^{(0)})_{,k} + \psi \Big(\rho V^2 u_{i,11}^{(0)} -2 \rho V
\fr{\prt^2 u_i^{(0)}}{\prt t \prt X} + \rho \fr{\prt^2 u_i^{(0)}}{\prt t^2}\Big) .
$$


\vspace{.5in}

\begin{thebibliography}{1}
\bibitem{MR}
Morrissey, J W and Rice, J R (1998)
{\it Crack front waves},
J. Mech. Phys. Solids {\bf 46}, 467--487.

\bibitem{RF}
Ramanathan, S and Fisher, D S (1997)
{\it Dynamics and instabilities of planar tensile cracks in heterogeneous media},
Phys. Rev. Lett. {\bf 79}, 877--880.

\bibitem{WM1}
Willis, J R  and Movchan, A B (1995)
{\it Dynamic weight functions for a moving
crack. I. Mode I loading}, J. Mech. Phys. Solids {\bf 43}, 319--341.


\bibitem{SCF}
Sharon, E, Cohen, G and Fineberg, J (2002)
{\it Propagating solitary waves along a rapidly moving
crack front}, Nature {\bf 410}, 68--71.

\bibitem{WM_visco}
Willis, J R and Movchan, A B (2001) {\it The influence of viscoelasticity
on crack front waves}, J. Mech. Phys. Solids, {\bf 49}, 2177-2189.

\bibitem{WMW}
Woolfries, S, Movchan, A B and Willis, J R (2002)
{\it Perturbation of a dynamic planar crack moving in a model viscoelastic solid},
Int. J. Solids Struct. {\bf 39}, 5409--5426.

\bibitem{OMW1}
Obrezanova, O, Movchan, A B and Willis, J R (2002)
{\it Dynamic stability of a propagating crack},
J. Mech. Phys. Solids,  {\bf 50}, 2637--2668.

\bibitem{LL}
Lazarus, V  and  Leblond, J-B (1998)
{\it Crack paths under mixed mode (I + III) or (I + II + III) loadings},
C. R. Acad. Sci. Paris, Series IIB, 
{\bf 326}, Issue 3, 171--177.


\bibitem{W_corr}
Willis, J R (2003) {\it Dynamic perturbation of a propagating crack:
implications for crack stability}, Asymptotics, Singularities and
Homogenization in Problems of Mechanics,
edited by A. B. Movchan, Kluwer, Dordrecht.

\bibitem{MW2}
Movchan, A B and Willis, J R (1995) {\it Dynamic weight functions for
a moving crack. II. Shear loading},  J. Mech. Phys. Solids {\bf 43},
1369--1383.

\bibitem{WM97}
Willis, J R  and Movchan, A B (1997) {\it Three-dimensional dynamic
perturbation of a propagating crack}, J. Mech. Phys. Solids {\bf 45},
591--610.

\bibitem{W}
Willis, J R (1999) {\it Asymptotic analysis in fracture: An update},
Int. J. Fract. {\bf 100}, 85--103.



\bibitem{F}
Freund, L B (1990) {\it Dynamic Fracture Mechanics.}
Cambridge: Cambridge University Press.

\bibitem{Oleaga}
Oleaga, G (2003) {\it On the dynamics of cracks in three dimensions},
J. Mech. Phys. Solids, {\bf 51}, 169-185.

\bibitem{B} Bueckner, H F (1987) {\it Weight functions and fundamental solutions
for the penny shaped and half-plane crack in three space},
Int. J. Solids Struct. {\bf 23}, 57--93.





\end{thebibliography}
\end{document}